\documentclass[11pt]{amsart}
\usepackage{amsmath,amssymb,stmaryrd}
\newtheorem{theorem}{Theorem}
\newtheorem{proposition}[theorem]{Proposition}
\newtheorem{corollary}[theorem]{Corollary}
\newtheorem{lemma}[theorem]{Lemma}
\begin{document}

\title[Ricci flow in higher dimensions]{A general convergence result for the Ricci flow in higher dimensions}
\author{Simon Brendle}
\address{Department of Mathematics \\
                 Stanford University \\
                 Stanford, CA 94305}
\thanks{This project was supported by the Alfred P. Sloan Foundation and by the National Science Foundation under grant DMS-0605223.}
\maketitle

\section{Introduction}

In this paper, we study the longterm behavior of the Ricci flow in higher dimensions. A one-parameter family of metrics 
$g(t)$ is a solution to the Ricci flow if 
\[\frac{\partial}{\partial t} g(t) = -2 \, \text{\rm Ric}_{g(t)},\] 
where $\text{\rm Ric}_{g(t)}$ denotes the Ricci tensor of $g(t)$ (cf. \cite{Hamilton1}). Moreover, $g(t)$ is a solution to the normalized Ricci flow if 
\[\frac{\partial}{\partial t} g(t) = -2 \, \text{\rm Ric}_{g(t)} + \frac{2}{n} \, r_{g(t)} \, g(t),\] 
where $r_{g(t)}$ denotes the mean value of the scalar curvature of $g(t)$. In a joint work with R.~Schoen, we proved the following theorem:

\begin{theorem}[\cite{Brendle-Schoen}, Theorem 3] 
\label{convergence.1}
Let $(M,g_0)$ be a compact Riemannian manifold of dimension $n \geq 4$. Assume that 
\begin{equation} 
\label{curvature.condition.1}
R_{1313} + \lambda^2 \, R_{1414} + \mu^2 \, R_{2323} + \lambda^2\mu^2 \, R_{2424} - 2\lambda\mu \, R_{1234} > 0 
\end{equation}
for all orthonormal four-frames $\{e_1,e_2,e_3,e_4\}$ and all $\lambda,\mu \in [-1,1]$. Then the normalized Ricci flow with initial metric $g_0$ exists for all time and converges to a constant curvature metric as $t \to \infty$.
\end{theorem}

An immediate consequence of Theorem \ref{convergence.1} is the Differentiable Sphere Theorem: if $(M,g_0)$ has strictly $1/4$-pinched sectional curvatures, then $M$ is diffeomorphic to a spherical space form. We refer to \cite{Brendle-Schoen} for a discussion of the history of this problem.

In this paper, we weaken the curvature assumption in Theorem \ref{convergence.1}. Our main result is the following:

\begin{theorem}
\label{convergence.2}
Let $(M,g_0)$ be a compact Riemannian manifold of dimension $n \geq 4$. Assume that 
\begin{equation} 
\label{curvature.condition.2}
R_{1313} + \lambda^2 \, R_{1414} + R_{2323} + \lambda^2 \, R_{2424} - 2\lambda \, R_{1234} > 0 
\end{equation}
for all orthonormal four-frames $\{e_1,e_2,e_3,e_4\}$ and all $\lambda \in [-1,1]$. Then the normalized Ricci flow with initial metric $g_0$ exists for all time and converges to a constant curvature metric as $t \to \infty$.
\end{theorem}

C.~B\"ohm and B.~Wilking \cite{Bohm-Wilking} have shown that the normalized Ricci flow deforms metrics with $2$-positive 
curvature operator to constant curvature metrics. It is easy to see that every manifold with $2$-positive curvature operator satisfies condition (\ref{curvature.condition.2}). Hence, the main theorem in \cite{Bohm-Wilking} is a subcase of Theorem \ref{convergence.2}. 

The conditions (\ref{curvature.condition.1}) and (\ref{curvature.condition.2}) are closely related to the notion of 
positive isotropic curvature. To explain this, suppose that $M$ is a Riemannian manifold of dimension $n \geq 4$. We say that $M$ has nonnegative isotropic curvature if 
\[R_{1313} + R_{1414} + R_{2323} + R_{2424} - 2 \, R_{1234} \geq 0\] 
for all orthonormal four-frames $\{e_1,e_2,e_3,e_4\}$ (cf. \cite{Micallef-Moore}, \cite{Micallef-Wang}). The product $M \times \mathbb{R}$ has nonnegative isotropic curvature if and only if 
\[R_{1313} + \lambda^2 \, R_{1414} + R_{2323} + \lambda^2 \, R_{2424} - 2\lambda \, R_{1234} \geq 0\] 
for all orthonormal four-frames $\{e_1,e_2,e_3,e_4\}$ and all $\lambda \in [-1,1]$ 
(see Proposition \ref{characterization.of.C.tilde} below). Similarly, the product $M \times \mathbb{R}^2$ has nonnegative isotropic curvature if and only if 
\[R_{1313} + \lambda^2 \, R_{1414} + \mu^2 \, R_{2323} + \lambda^2\mu^2 \, R_{2424} - 2\lambda\mu \, R_{1234} \geq 0\] 
for all orthonormal four-frames $\{e_1,e_2,e_3,e_4\}$ and all $\lambda,\mu \in [-1,1]$ (see \cite{Brendle-Schoen}, Proposition 21). 

The curvature conditions (\ref{curvature.condition.1}) and (\ref{curvature.condition.2}) are void in dimension less than $4$. However, the condition that $M \times \mathbb{R}$ has nonnegative isotropic curvature makes sense for all $n \geq 3$, and the condition that $M \times \mathbb{R}^2$ has nonnegative isotropic curvature makes sense for all $n \geq 2$. A three-manifold $M$ has nonnegative Ricci curvature if and only if $M \times \mathbb{R}$ has nonnegative isotropic curvature. Moreover, a three-manifold $M$ has nonnegative sectional curvature if and only if $M \times \mathbb{R}^2$ has nonnegative isotropic curvature. Thus, Theorem \ref{convergence.2} can be viewed as a generalization of a theorem of R.~Hamilton on three-manifolds with positive Ricci curvature (see \cite{Hamilton1}). Combining the two results, we obtain:

\begin{theorem}
\label{M.cross.R.PIC}
Let $(M,g_0)$ be a compact Riemannian manifold of dimension $n \geq 3$. If $(M,g_0) \times \mathbb{R}$ has positive isotropic curvature, then the normalized Ricci flow with initial metric $g_0$ exists for all time and converges to a constant curvature metric as $t \to \infty$.
\end{theorem}

R.~Hamilton \cite{Hamilton3} has shown that the Ricci flow preserves positive isotropic curvature in dimension $4$. Moreover, Hamilton proved that, in dimension $4$, a solution to the Ricci flow with positive isotropic curvature develops only "neck-like" singularities. More recently, it was shown that positive isotropic curvature is preserved by the Ricci flow in all dimensions. This result was proved independently in \cite{Brendle-Schoen} and \cite{Nguyen}. It is an open question whether the analysis of singularities in \cite{Hamilton3} carries over to higher dimensions. We hope that Theorem \ref{M.cross.R.PIC} will shed light on this question.

In Section 2, we consider the condition that $M \times \mathbb{R}$ has nonnegative isotropic curvature. This condition defines a convex cone $\tilde{C}$ in the space of algebraic curvature operators, which is preserved by the Hamilton ODE. 

In Section 3, we consider the condition that $M \times S^2(1)$ has nonnegative isotropic curvature. This defines a convex set $E$ in the space of algebraic curvature operators. It is easy to see that $\hat{C} \subset E \subset \tilde{C}$, where $\hat{C}$ denotes the cone introduced in \cite{Brendle-Schoen}. Using results from \cite{Brendle-Schoen}, we show that the set $E$ is invariant under the Hamilton ODE. This fact is the main ingredient in the proof of Theorem 
\ref{convergence.2}. 

In Section 4, we complete the proof of Theorem \ref{convergence.2} by constructing a suitable pinching set for the Hamilton ODE.

\section{The cone $\tilde{C}$}

Let $R$ be an algebraic curvature operator on $\mathbb{R}^n$. We define an algebraic curvature operator $\tilde{R}$ on $\mathbb{R}^n \times \mathbb{R}$ by 
\[\tilde{R}(\tilde{v}_1,\tilde{v}_2,\tilde{v}_3,\tilde{v}_4) = R(v_1,v_2,v_3,v_4)\] 
for all vectors $\tilde{v}_j = (v_j,x_j) \in \mathbb{R}^n \times \mathbb{R}$. We denote by $\tilde{C}$ the set of all algebraic curvature operators on $\mathbb{R}^n$ with the property that $\tilde{R}$ has nonnegative isotropic curvature:
\[\tilde{C} = \{R \in S_B^2(\mathfrak{so}(n)): \text{$\tilde{R}$ has nonnegative isotropic curvature}\}.\] 
Clearly, $\tilde{C}$ is closed, convex, and $O(n)$-invariant. Moreover, it follows from the results in \cite{Brendle-Schoen} that $\tilde{C}$ is invariant under the Hamilton ODE $\frac{d}{dt} R = Q(R)$. The cone $\tilde{C}$ can be characterized as follows:

\begin{proposition} 
\label{characterization.of.C.tilde}
Let $R$ be an algebraic curvature operator on $\mathbb{R}^n$, and let $\tilde{R}$ be the induced curvature operator on $\mathbb{R}^n \times \mathbb{R}$. The following statements are equivalent: 
\begin{itemize}
\item[(i)] $\tilde{R}$ has nonnegative isotropic curvature.
\item[(ii)] For all orthonormal four-frames $\{e_1,e_2,e_3,e_4\}$ and all $\lambda \in [-1,1]$, we have 
\begin{align*} 
&R(e_1,e_3,e_1,e_3) + \lambda^2 \, R(e_1,e_4,e_1,e_4) \\ 
&+ R(e_2,e_3,e_2,e_3) + \lambda^2 \, R(e_2,e_4,e_2,e_4) - 2\lambda \, R(e_1,e_2,e_3,e_4) \geq 0. 
\end{align*} 
\end{itemize}
\end{proposition}

\textbf{Proof.} 
Assume first that $\tilde{R}$ has nonnegative isotropic curvature. Let 
$\{e_1,e_2,e_3,e_4\}$ be an orthonormal four-frame in $\mathbb{R}^n$, and let $\lambda \in [-1,1]$. We define 
\[\begin{array}{l@{\qquad\qquad}l} 
\tilde{e}_1 = (e_1,0) & \tilde{e}_2 = (e_2,0) \\ 
\tilde{e}_3 = (e_3,0) & \tilde{e}_4 = (\lambda e_4,\sqrt{1-\lambda^2}). 
\end{array}\] 
Since $\tilde{R}$ has nonnegative isotropic curvature, we have 
\begin{align*} 
0 &\leq \tilde{R}(\tilde{e}_1,\tilde{e}_3,\tilde{e}_1,\tilde{e}_3) + \tilde{R}(\tilde{e}_1,\tilde{e}_4,\tilde{e}_1,\tilde{e}_4) \\ 
&+ \tilde{R}(\tilde{e}_2,\tilde{e}_3,\tilde{e}_2,\tilde{e}_3) + \tilde{R}(\tilde{e}_2,\tilde{e}_4,\tilde{e}_2,\tilde{e}_4) - 2 \, \tilde{R}(\tilde{e}_1,\tilde{e}_2,\tilde{e}_3,\tilde{e}_4) \\ 
&= R(e_1,e_3,e_1,e_3) + \lambda^2 \, R(e_1,e_4,e_1,e_4) \\ 
&+ R(e_2,e_3,e_2,e_3) + \lambda^2 \, R(e_2,e_4,e_2,e_4) - 2\lambda \, R(e_1,e_2,e_3,e_4), 
\end{align*} 
as claimed.

Conversely, assume that (ii) holds. We claim that $\tilde{R}$ has nonnegative isotropic curvature. 
Let $\{\tilde{e}_1,\tilde{e}_2,\tilde{e}_3,\tilde{e}_4\}$ be an orthonormal four-frame in $\mathbb{R}^n 
\times \mathbb{R}$. We write $\tilde{e}_j = (v_j,x_j)$, where $v_j \in \mathbb{R}^n$ and 
$x_j \in \mathbb{R}$. Moreover, we define 
\begin{align*} 
\varphi &= v_1 \wedge v_3 + v_4 \wedge v_2 \\ 
\psi &= v_1 \wedge v_4 + v_2 \wedge v_3. 
\end{align*}
Clearly, $\varphi \wedge \varphi = \psi \wedge \psi$ and $\varphi \wedge \psi = 0$. 
Using the relation $\langle v_i,v_j \rangle + x_ix_j = \delta_{ij}$, we obtain 
\begin{align*} 
|\varphi|^2 - |\psi|^2 
&= |v_1 \wedge v_3|^2 + |v_4 \wedge v_2|^2 - |v_1 \wedge v_4|^2 - |v_2 \wedge v_3|^2 \\ 
&+ 2 \, \langle v_1 \wedge v_3,v_4 \wedge v_2 \rangle - 2 \, \langle v_1 \wedge v_4,v_2 \wedge v_3 \rangle \\ 
&= (|v_1|^2 - |v_2|^2)(|v_3|^2 - |v_4|^2) - 4 \, \langle v_1,v_2 \rangle \, \langle v_3,v_4 \rangle \\ 
&- (\langle v_1,v_3 \rangle - \langle v_2,v_4 \rangle)^2 + (\langle v_1,v_4 \rangle^2 + \langle v_2,v_3 \rangle)^2 \\ 
&= (x_1^2 - x_2^2)(x_3^2 - x_4^2) - 4x_1x_2x_3x_4 \\ 
&- (x_1x_3 - x_2x_4)^2 + (x_1x_4 + x_2x_3)^2 \\ 
&= 0
\end{align*} 
and 
\begin{align*} 
\langle \varphi,\psi \rangle 
&= \langle v_1 \wedge v_3,v_1 \wedge v_4 \rangle + \langle v_1 \wedge v_3,v_2 \wedge v_3 \rangle \\ 
&+ \langle v_4 \wedge v_2,v_1 \wedge v_4 \rangle + \langle v_4 \wedge v_2,v_2 \wedge v_3 \rangle \\ 
&= (|v_1|^2 - |v_2|^2) \, \langle v_3,v_4 \rangle + (|v_3|^2 - |v_4|^2) \, \langle v_1,v_2 \rangle \\ 
&- (\langle v_1,v_3 \rangle - \langle v_2,v_4 \rangle) \, (\langle v_1,v_4 \rangle + \langle v_2,v_3 \rangle) \\ 
&= (x_1^2 - x_2^2) \, x_3x_4 + (x_3^2 - x_4^2) \, x_1x_2 \\ 
&- (x_1x_3 - x_2x_4) \, (x_1x_4 + x_2x_3) \\ 
&= 0. 
\end{align*} 
By Lemma 19 in \cite{Brendle-Schoen}, we can find an orthonormal four-frame $\{e_1,e_2,e_3,e_4\}$ in $\mathbb{R}^n$ and real numbers $a_1,a_2,b_1,b_2$ such that $a_1^2 + a_2^2 = b_1^2 + b_2^2$, $a_1a_2 = b_1b_2$, and 
\begin{align*} 
\varphi &= a_1 \, e_1 \wedge e_3 + a_2 \, e_4 \wedge e_2 \\ 
\psi &= b_1 \, e_1 \wedge e_4 + b_2 \, e_2 \wedge e_3. 
\end{align*} 
Clearly, $(a_1^2 - b_1^2)(a_1^2 - b_2^2) = 0$. Without loss of generality, we may assume that $a_1^2 = b_2^2$. (Otherwise, we replace $\{e_1,e_2,e_3,e_4\}$ by $\{e_3,e_4,e_1,e_2\}$.) This implies $a_2^2 = b_1^2$. Using the first Bianchi identity, we obtain 
\begin{align*} 
R(\varphi,\varphi) + R(\psi,\psi) &= a_1^2 \, R(e_1,e_3,e_1,e_3) + a_2^2 \, R(e_1,e_4,e_1,e_4) \\ 
&+ a_1^2 \, R(e_2,e_3,e_2,e_3) + a_2^2 \, R(e_2,e_4,e_2,e_4) \\ 
&- 2a_1a_2 \, R(e_1,e_2,e_3,e_4). 
\end{align*} 
The condition (ii) implies that the right hand side is nonnegative. Thus, we conclude that 
\begin{align*} 
&\tilde{R}(\tilde{e}_1,\tilde{e}_3,\tilde{e}_1,\tilde{e}_3) + \tilde{R}(\tilde{e}_1,\tilde{e}_4,\tilde{e}_1,\tilde{e}_4) \\ 
&+ \tilde{R}(\tilde{e}_2,\tilde{e}_3,\tilde{e}_2,\tilde{e}_3) + \tilde{R}(\tilde{e}_2,\tilde{e}_4,\tilde{e}_2,\tilde{e}_4) - 2 \, \tilde{R}(\tilde{e}_1,\tilde{e}_2,\tilde{e}_3,\tilde{e}_4) \\ 
&= R(\varphi,\varphi) + R(\psi,\psi) \geq 0. 
\end{align*} 
Hence, $\tilde{R}$ has nonnegative isotropic curvature. \\

\section{A new invariant curvature condition}

Let $R$ be an algebraic curvature operator on $\mathbb{R}^n$. Following Hamilton \cite{Hamilton2}, we define  
\[Q(R)_{ijkl} = R_{ijpq} \, R_{klpq} + 2 \, R_{ipkq} \, R_{jplq} - 2 \, R_{iplq} \, R_{jpkq}.\] 
It is straightforward to verify that $Q(R)$ is an algebraic curvature tensor. As in \cite{Hamilton2}, we write 
$Q(R) = R^2 + R^\#$, where $R^2$ and $R^\#$ are defined by 
\begin{align*} 
(R^2)_{ijkl} &= R_{ijpq} \, R_{klpq} \\ 
(R^\#)_{ijkl} &= 2 \, R_{ipkq} \, R_{jplq} - 2 \, R_{iplq} \, R_{jpkq}. 
\end{align*} 
Note that $R^2$ and $R^\#$ do not satisfy the first Bianchi identity, but 
$R^2 + R^\#$ does. The following lemma is a consequence of Corollary 10 in 
\cite{Brendle-Schoen}, and plays a key role in our analysis:

\begin{lemma} 
\label{R.sharp.term.1}
Let $R$ be an algebraic curvature operator on $\mathbb{R}^n$ with nonnegative isotropic curvature. Moreover, suppose that 
$\{e_1,e_2,e_3,e_4\}$ is an orthonormal four-frame in $\mathbb{R}^n$ satisfying 
\begin{align} 
&R(e_1,e_3,e_1,e_3) + R(e_1,e_4,e_1,e_4) \notag \\ 
&+ R(e_2,e_3,e_2,e_3) + R(e_2,e_4,e_2,e_4) - 2 \, R(e_1,e_2,e_3,e_4) = 0. 
\end{align}
Then 
\begin{align} 
&R^\#(e_1,e_3,e_1,e_3) + R^\#(e_1,e_4,e_1,e_4) \notag \\ 
&+ R^\#(e_2,e_3,e_2,e_3) + R^\#(e_2,e_4,e_2,e_4) \\ 
&+ 2 \, R^\#(e_1,e_3,e_4,e_2) + 2 \, R^\#(e_1,e_4,e_2,e_3) \geq 0. \notag 
\end{align}
\end{lemma}

\textbf{Proof.} 
We extend $\{e_1,e_2,e_3,e_4\}$ to an orthonormal basis $\{e_1, \hdots, e_n\}$ of $\mathbb{R}^n$. 
Using the first Bianchi identity, we obtain 
\begin{align*} 
&R^\#(e_1,e_3,e_4,e_2) + R^\#(e_1,e_4,e_2,e_3) \\ 
&= 2 \, R(e_1,e_p,e_4,e_q) \, R(e_3,e_p,e_2,e_q) - 2 \, R(e_1,e_p,e_3,e_q) \, R(e_4,e_p,e_2,e_q) \\ 
&+ 2 \, R(e_1,e_p,e_2,e_q) \, R(e_4,e_p,e_3,e_q) - 2 \, R(e_1,e_p,e_2,e_q) \, R(e_3,e_p,e_4,e_q) \\ 
&= 2 \, R(e_1,e_p,e_4,e_q) \, R(e_3,e_p,e_2,e_q) - 2 \, R(e_1,e_p,e_3,e_q) \, R(e_4,e_p,e_2,e_q) \\ 
&- R(e_1,e_2,e_p,e_q) \, R(e_3,e_4,e_p,e_q). 
\end{align*} 
This implies 
\begin{align*} 
&R^\#(e_1,e_3,e_1,e_3) + R^\#(e_1,e_4,e_1,e_4) \\ 
&+ R^\#(e_2,e_3,e_2,e_3) + R^\#(e_2,e_4,e_2,e_4) \\ 
&+ 2 \, R^\#(e_1,e_3,e_4,e_2) + 2 \, R^\#(e_1,e_4,e_2,e_3) \\ 
&= 2 \, R(e_1,e_p,e_1,e_q) \, R(e_3,e_p,e_3,e_q) - 2 \, R(e_1,e_p,e_3,e_q) \, R(e_3,e_p,e_1,e_q) \\ 
&+ 2 \, R(e_1,e_p,e_1,e_q) \, R(e_4,e_p,e_4,e_q) - 2 \, R(e_1,e_p,e_4,e_q) \, R(e_4,e_p,e_1,e_q) \\ 
&+ 2 \, R(e_2,e_p,e_2,e_q) \, R(e_3,e_p,e_3,e_q) - 2 \, R(e_2,e_p,e_3,e_q) \, R(e_3,e_p,e_2,e_q) \\ 
&+ 2 \, R(e_2,e_p,e_2,e_q) \, R(e_4,e_p,e_4,e_q) - 2 \, R(e_2,e_p,e_4,e_q) \, R(e_4,e_p,e_2,e_q) \\ 
&+ 4 \, R(e_1,e_p,e_4,e_q) \, R(e_3,e_p,e_2,e_q) - 4 \, R(e_1,e_p,e_3,e_q) \, R(e_4,e_p,e_2,e_q) \\ 
&- 2 \, R(e_1,e_2,e_p,e_q) \, R(e_3,e_4,e_p,e_q). 
\end{align*} 
Rearranging terms yields 
\begin{align*} 
&R^\#(e_1,e_3,e_1,e_3) + R^\#(e_1,e_4,e_1,e_4) \\ 
&+ R^\#(e_2,e_3,e_2,e_3) + R^\#(e_2,e_4,e_2,e_4) \\ 
&+ 2 \, R^\#(e_1,e_3,e_4,e_2) + 2 \, R^\#(e_1,e_4,e_2,e_3) \\ 
&= 2 \, (R(e_1,e_p,e_1,e_q) + R(e_2,e_p,e_2,e_q)) \, (R(e_3,e_p,e_3,e_q) + R(e_4,e_p,e_4,e_q)) \\ 
&- 2 \, R(e_1,e_2,e_p,e_q) \, R(e_3,e_4,e_p,e_q) \\ 
&- 2 \, (R(e_1,e_p,e_3,e_q) + R(e_2,e_p,e_4,e_q)) \, (R(e_3,e_p,e_1,e_q) + R(e_4,e_p,e_2,e_q)) \\ 
&+ 2 \, (R(e_1,e_p,e_4,e_q) - R(e_2,e_p,e_3,e_q)) \, (R(e_4,e_p,e_1,e_q) - R(e_3,e_p,e_2,e_q)), 
\end{align*}
and the right hand side is nonnegative by Corollary 10 in \cite{Brendle-Schoen}. \\

Given any algebraic curvature operator $R$ on $\mathbb{R}^n$, we define an algebraic curvature operator $S$ on $\mathbb{R}^n \times \mathbb{R}^2$ by 
\[S(\hat{v}_1,\hat{v}_2,\hat{v}_3,\hat{v}_4) = R(v_1,v_2,v_3,v_4) + \langle x_1,x_3 \rangle \, \langle x_2,x_4 \rangle - \langle x_1,x_4 \rangle \, \langle x_2,x_3 \rangle\] 
for all vectors $\hat{v}_j = (v_j,x_j) \in \mathbb{R}^n \times \mathbb{R}^2$. A straightforward calculation yields:

\begin{lemma}
\label{R.sharp.term.2}
Let $R$ be an algebraic curvature operator on $\mathbb{R}^n$, and let $S$ be the induced curvature operator on $\mathbb{R}^n \times \mathbb{R}^2$. Then 
\[S^\#(\hat{v}_1,\hat{v}_2,\hat{v}_3,\hat{v}_4) = R^\#(v_1,v_2,v_3,v_4)\] 
for all vectors $\hat{v}_j = (v_j,x_j) \in \mathbb{R}^n \times \mathbb{R}^2$. 
\end{lemma}

Let $E$ be the set of all algebraic curvature operators on $\mathbb{R}^n$ with the property that the induced curvature operator $S$ on $\mathbb{R}^n \times \mathbb{R}^2$ has nonnegative isotropic curvature:
\[E = \{R \in S_B^2(\mathfrak{so}(n)): \text{$S$ has nonnegative isotropic curvature}\}\] 
It is easy to see that $E$ is closed, convex, and $O(n)$-invariant.

\begin{proposition} 
\label{algebraic.fact}
Let $R$ be an algebraic curvature operator on $\mathbb{R}^n$, and let $S$ be the induced curvature operator on $\mathbb{R}^n \times \mathbb{R}^2$. The following statements are equivalent: 
\begin{itemize}
\item[(i)] $S$ has nonnegative isotropic curvature.
\item[(ii)] For all orthonormal four-frames $\{e_1,e_2,e_3,e_4\}$ and all $\lambda,\mu \in [-1,1]$, we have 
\begin{align*} 
&R(e_1,e_3,e_1,e_3) + \lambda^2 \, R(e_1,e_4,e_1,e_4) \\ 
&+ \mu^2 \, R(e_2,e_3,e_2,e_3) + \lambda^2\mu^2 \, R(e_2,e_4,e_2,e_4) \\ 
&- 2\lambda\mu \, R(e_1,e_2,e_3,e_4) + (1 - \lambda^2) \, (1 - \mu^2) \geq 0. 
\end{align*} 
\end{itemize}
\end{proposition}

\textbf{Proof.} 
Assume first that $S$ has nonnegative isotropic curvature. Let 
$\{e_1,e_2,e_3,e_4\}$ be an orthonormal four-frame in $\mathbb{R}^n$, and let $\lambda,\mu \in [-1,1]$. We define 
\[\begin{array}{l@{\qquad\qquad}l} 
\hat{e}_1 = (e_1,0,0) & \hat{e}_2 = (\mu e_2,0,\sqrt{1-\mu^2}) \\ 
\hat{e}_3 = (e_3,0,0) & \hat{e}_4 = (\lambda e_4,\sqrt{1-\lambda^2},0). 
\end{array}\] 
Clearly, the vectors $\{\hat{e}_1,\hat{e}_2,\hat{e}_3,\hat{e}_4\}$ form an orthonormal four-frame in 
$\mathbb{R}^n \times \mathbb{R}^2$. Since $S$ has nonnegative isotropic curvature, we have 
\begin{align*} 
0 &\leq S(\hat{e}_1,\hat{e}_3,\hat{e}_1,\hat{e}_3) + S(\hat{e}_1,\hat{e}_4,\hat{e}_1,\hat{e}_4) \\ 
&+ S(\hat{e}_2,\hat{e}_3,\hat{e}_2,\hat{e}_3) + S(\hat{e}_2,\hat{e}_4,\hat{e}_2,\hat{e}_4) - 2 \, S(\hat{e}_1,\hat{e}_2,\hat{e}_3,\hat{e}_4) \\ 
&= R(e_1,e_3,e_1,e_3) + \lambda^2 \, R(e_1,e_4,e_1,e_4) \\ 
&+ \mu^2 \, R(e_2,e_3,e_2,e_3) + \lambda^2\mu^2 \, R(e_2,e_4,e_2,e_4) \\ 
&- 2\lambda\mu \, R(e_1,e_2,e_3,e_4) + (1 - \lambda^2) \, (1 - \mu^2), 
\end{align*} 
as claimed.

Conversely, assume that (ii) holds. We claim that $S$ has nonnegative isotropic curvature. 
Let $\{\hat{e}_1,\hat{e}_2,\hat{e}_3,\hat{e}_4\}$ be an orthonormal four-frame in $\mathbb{R}^n 
\times \mathbb{R}^2$. We write $\hat{e}_j = (v_j,x_j)$, where $v_j \in \mathbb{R}^n$ and $x_j \in \mathbb{R}^2$. 
Let $V$ be a four-dimensional subspace of $\mathbb{R}^n$ containing $\{v_1,v_2,v_3,v_4\}$. We define 
\begin{align*} 
\varphi &= v_1 \wedge v_3 + v_4 \wedge v_2 \in \wedge^2 V \\ 
\psi &= v_1 \wedge v_4 + v_2 \wedge v_3 \in \wedge^2 V. 
\end{align*}
Clearly, $\varphi \wedge \varphi = \psi \wedge \psi$ and $\varphi \wedge \psi = 0$. 
By Lemma 20 in \cite{Brendle-Schoen}, there exist an orthonormal basis $\{e_1,e_2,e_3,e_4\}$ of $V$ and real numbers $a_1,a_2,b_1,b_2,\theta$ such that $a_1a_2 = b_1b_2$ and 
\begin{align*} 
\tilde{\varphi} &:= \cos \theta \, \varphi + \sin \theta \, \psi = a_1 \, e_1 \wedge e_3 + a_2 \, e_4 \wedge e_2 \\ 
\tilde{\psi} &:= -\sin \theta \, \varphi + \cos \theta \, \psi = b_1 \, e_1 \wedge e_4 + b_2 \, e_2 \wedge e_3. 
\end{align*} 
This implies
\begin{align*} 
R(\varphi,\varphi) + R(\psi,\psi) &= R(\tilde{\varphi},\tilde{\varphi}) + R(\tilde{\psi},\tilde{\psi}) \\ 
&= a_1^2 \, R(e_1,e_3,e_1,e_3) + b_1^2 \, R(e_1,e_4,e_1,e_4) \\ 
&+ b_2^2 \, R(e_2,e_3,e_2,e_3) + a_2^2 \, R(e_2,e_4,e_2,e_4) \\ 
&- 2a_1a_2 \, R(e_1,e_2,e_3,e_4). 
\end{align*}
Using the identity $\langle v_i,v_j \rangle + \langle x_i,x_j \rangle = \delta_{ij}$, we obtain 
\begin{align*} 
|\varphi|^2 - |\psi|^2 
&= (|v_1|^2 - |v_2|^2)(|v_3|^2 - |v_4|^2) - 4 \, \langle v_1,v_2 \rangle \, \langle v_3,v_4 \rangle \\ 
&- (\langle v_1,v_3 \rangle - \langle v_2,v_4 \rangle)^2 + (\langle v_1,v_4 \rangle^2 + \langle v_2,v_3 \rangle)^2 \\ 
&= (|x_1|^2 - |x_2|^2)(|x_3|^2 - |x_4|^2) - 4 \, \langle x_1,x_2 \rangle \, \langle x_3,x_4 \rangle \\ 
&- (\langle x_1,x_3 \rangle - \langle x_2,x_4 \rangle)^2 + (\langle x_1,x_4 \rangle + \langle x_2,x_3 \rangle)^2 \\ 
&= |x_1 \wedge x_3 + x_4 \wedge x_2|^2 - |x_1 \wedge x_4 + x_2 \wedge x_3|^2 
\end{align*} 
and 
\begin{align*} 
\langle \varphi,\psi \rangle 
&= (|v_1|^2 - |v_2|^2) \, \langle v_3,v_4 \rangle + (|v_3|^2 - |v_4|^2) \, \langle v_1,v_2 \rangle \\ 
&- (\langle v_1,v_3 \rangle - \langle v_2,v_4 \rangle) \, (\langle v_1,v_4 \rangle + \langle v_2,v_3 \rangle) \\ 
&= (|x_1|^2 - |x_2|^2) \, \langle x_3,x_4 \rangle + (|x_3|^2 - |x_4^2|) \, \langle x_1,x_2 \rangle \\ 
&- (\langle x_1,x_3 \rangle - \langle x_2,x_4 \rangle) \, (\langle x_1,x_4 \rangle + \langle x_2,x_3 \rangle) \\ 
&= \langle x_1 \wedge x_3 + x_4 \wedge x_2,x_1 \wedge x_4 + x_2 \wedge x_3 \rangle. 
\end{align*} From this we deduce that 
\begin{align*} 
&(|x_1 \wedge x_3 + x_4 \wedge x_2|^2 + |x_1 \wedge x_4 + x_2 \wedge x_3|^2)^2 \\ 
&= (|x_1 \wedge x_3 + x_4 \wedge x_2|^2 - |x_1 \wedge x_4 + x_2 \wedge x_3|^2)^2 \\ 
&+ 4 \, |x_1 \wedge x_3 + x_4 \wedge x_2|^2 \, |x_1 \wedge x_4 + x_2 \wedge x_3|^2 \\ 
&\geq (|x_1 \wedge x_3 + x_4 \wedge x_2|^2 - |x_1 \wedge x_4 + x_2 \wedge x_3|^2)^2 \\ 
&+ 4 \, \langle x_1 \wedge x_3 + x_4 \wedge x_2,x_1 \wedge x_4 + x_2 \wedge x_3 \rangle^2 \\ 
&= (|\varphi|^2 - |\psi|^2)^2 + 4 \, \langle \varphi,\psi \rangle^2 \\ 
&= (|\tilde{\varphi}|^2 - |\tilde{\psi}|^2)^2 + 4 \, \langle \tilde{\varphi},\tilde{\psi} \rangle^2 \\ 
&= (a_1^2 + a_2^2 - b_1^2 - b_2^2)^2. 
\end{align*} 
Putting these facts together, we obtain 
\begin{align*} 
&R(\varphi,\varphi) + R(\psi,\psi) + |x_1 \wedge x_3 + x_4 \wedge x_2|^2 + |x_1 \wedge x_4 + x_2 \wedge x_3|^2 \\ 
&\geq a_1^2 \, R(e_1,e_3,e_1,e_3) + b_1^2 \, R(e_1,e_4,e_1,e_4) \\ 
&+ b_2^2 \, R(e_2,e_3,e_2,e_3) + a_2^2 \, R(e_2,e_4,e_2,e_4) \\ 
&- 2a_1a_2 \, R(e_1,e_2,e_3,e_4) + |a_1^2 + a_2^2 - b_1^2 - b_2^2|. 
\end{align*} 
The condition (ii) implies that the right hand side is nonnegative. Thus, we conclude that 
\begin{align*} 
&S(\hat{e}_1,\hat{e}_3,\hat{e}_1,\hat{e}_3) + S(\hat{e}_1,\hat{e}_4,\hat{e}_1,\hat{e}_4) \\ 
&+ S(\hat{e}_2,\hat{e}_3,\hat{e}_2,\hat{e}_3) + S(\hat{e}_2,\hat{e}_4,\hat{e}_2,\hat{e}_4) 
- 2 \, S(\hat{e}_1,\hat{e}_2,\hat{e}_3,\hat{e}_4) \\ 
&= R(\varphi,\varphi) + R(\psi,\psi) + |x_1 \wedge x_3 + x_4 \wedge x_2|^2 + |x_1 \wedge x_4 + x_2 \wedge x_3|^2 \geq 0, 
\end{align*}
as claimed. \\

We next consider the cone $\hat{C}$ introduced in \cite{Brendle-Schoen}. Moreover, we denote by $I$ the 
curvature operator of the standard sphere, i.e. $I_{ijkl} = \delta_{ik} \, \delta_{jl} - \delta_{il} \, \delta_{jk}$. 
Using Proposition 21 in \cite{Brendle-Schoen}, we obtain:

\begin{corollary}
\label{inclusions}
If $R \in E$, then $R \in \tilde{C}$ and $R + I \in \hat{C}$. Moreover, we have $E + \hat{C} = E$.
\end{corollary}

We claim that the set $E$ is invariant under the ODE $\frac{d}{dt} R = Q(R)$. This is a consequence of the following algebraic fact: 

\begin{proposition} 
\label{invariance.1}
Let $R \in E$ be an algebraic curvature operator on $\mathbb{R}^n$. Moreover, let $\{e_1,e_2,e_3,e_4\}$ be an orthonormal four-frame in $\mathbb{R}^n$, and let $\lambda,\mu \in [-1,1]$. If 
\begin{align} 
\label{assumption}
&R(e_1,e_3,e_1,e_3) + \lambda^2 \, R(e_1,e_4,e_1,e_4) \notag \\ 
&+ \mu^2 \, R(e_2,e_3,e_2,e_3) + \lambda^2\mu^2 \, R(e_2,e_4,e_2,e_4) \\ 
&- 2\lambda\mu \, R(e_1,e_2,e_3,e_4) + (1 - \lambda^2) \, (1 - \mu^2) = 0, \notag 
\end{align} 
then we have 
\begin{align} 
&Q(R)(e_1,e_3,e_1,e_3) + \lambda^2 \, Q(R)(e_1,e_4,e_1,e_4) \notag \\ 
&+ \mu^2 \, Q(R)(e_2,e_3,e_2,e_3) + \lambda^2\mu^2 \, Q(R)(e_2,e_4,e_2,e_4) \\ 
&- 2\lambda\mu \, Q(R)(e_1,e_2,e_3,e_4) \geq 0. \notag 
\end{align}
\end{proposition}

\textbf{Proof.} 
Let $S$ be the curvature operator on $\mathbb{R}^n \times \mathbb{R}^2$ associated with $R$. We define an orthonormal four-frame $\{\hat{e}_1,\hat{e}_2,\hat{e}_3,\hat{e}_4\}$ in $\mathbb{R}^n \times \mathbb{R}^2$ by 
\[\begin{array}{l@{\qquad\qquad}l} 
\hat{e}_1 = (e_1,0,0) & \hat{e}_2 = (\mu e_2,0,\sqrt{1-\mu^2}) \\ 
\hat{e}_3 = (e_3,0,0) & \hat{e}_4 = (\lambda e_4,\sqrt{1-\lambda^2},0). 
\end{array}\] 
By assumption, $S$ has nonnegative isotropic curvature. Moreover, it follows from (\ref{assumption}) that 
\begin{align*} 
&S(\hat{e}_1,\hat{e}_3,\hat{e}_1,\hat{e}_3) + S(\hat{e}_1,\hat{e}_4,\hat{e}_1,\hat{e}_4) \\ 
&+ S(\hat{e}_2,\hat{e}_3,\hat{e}_2,\hat{e}_3) + S(\hat{e}_2,\hat{e}_4,\hat{e}_2,\hat{e}_4) - 2 \, S(\hat{e}_1,\hat{e}_2,\hat{e}_3,\hat{e}_4) = 0. 
\end{align*}
Hence, Lemma \ref{R.sharp.term.1} implies that 
\begin{align*} 
&S^\#(\hat{e}_1,\hat{e}_3,\hat{e}_1,e_3) + S^\#(\hat{e}_1,\hat{e}_4,\hat{e}_1,\hat{e}_4) \\ 
&+ S^\#(\hat{e}_2,\hat{e}_3,\hat{e}_2,\hat{e}_3) + S^\#(\hat{e}_2,\hat{e}_4,\hat{e}_2,\hat{e}_4) \\ 
&+ 2 \, S^\#(\hat{e}_1,\hat{e}_3,\hat{e}_4,\hat{e}_2) + 2 \, S^\#(\hat{e}_1,\hat{e}_4,\hat{e}_2,\hat{e}_3) \geq 0. 
\end{align*} 
Using Lemma \ref{R.sharp.term.2}, we obtain 
\begin{align} 
\label{estimate.1}
&R^\#(e_1,e_3,e_1,e_3) + \lambda^2 \, R^\#(e_1,e_4,e_1,e_4) \notag \\ 
&+ \mu^2 \, R^\#(e_2,e_3,e_2,e_3) + \lambda^2\mu^2 \, R^\#(e_2,e_4,e_2,e_4) \\ 
&+ 2\lambda\mu \, R^\#(e_1,e_3,e_4,e_2) + 2\lambda\mu \, R^\#(e_1,e_4,e_2,e_3) \geq 0. \notag 
\end{align} 
Moreover, we have 
\begin{align} 
\label{estimate.2}
&R^2(e_1,e_3,e_1,e_3) + \lambda^2 \, R^2(e_1,e_4,e_1,e_4) \notag \\ 
&+ \mu^2 \, R^2(e_2,e_3,e_2,e_3) + \lambda^2 \mu^2 \, R^2(e_2,e_4,e_2,e_4) \notag \\ 
&+ 2\lambda\mu \, R^2(e_1,e_3,e_4,e_2) + 2\lambda\mu \, R^2(e_1,e_4,e_2,e_3) \\ 
&= \sum_{p,q=1}^n \big [ R(e_1,e_3,e_p,e_q) - \lambda\mu \, R(e_2,e_4,e_p,e_q) \big ]^2 \notag \\ 
&+ \sum_{p,q=1}^n \big [ \lambda \, R(e_1,e_4,e_p,e_q) + \mu \, R(e_2,e_3,e_p,e_q) \big ]^2 \geq 0. \notag 
\end{align} 
Adding (\ref{estimate.1}) and (\ref{estimate.2}), we conclude that 
\begin{align} 
&Q(R)(e_1,e_3,e_1,e_3) + \lambda^2 \, Q(R)(e_1,e_4,e_1,e_4) \notag \\ 
&+ \mu^2 \, Q(R)(e_2,e_3,e_2,e_3) + \lambda^2\mu^2 \, Q(R)(e_2,e_4,e_2,e_4) \\ 
&+ 2\lambda\mu \, Q(R)(e_1,e_3,e_4,e_2) + 2\lambda\mu \, Q(R)(e_1,e_4,e_2,e_3) \geq 0. \notag 
\end{align} 
Since 
\[Q(R)(e_1,e_2,e_3,e_4) + Q(R)(e_1,e_3,e_4,e_2) + Q(R)(e_1,e_4,e_2,e_3) = 0,\] 
the assertion follows. \\

\begin{proposition} 
\label{invariance.2}
Suppose that $R(t)$, $t \in [0,T)$, is a solution of the ODE $\frac{d}{dt} R(t) = Q(R(t))$ with $R(0) \in E$. Then $R(t) \in E$ for all $t \in [0,T)$.
\end{proposition}

\textbf{Proof.} 
Fix $\varepsilon > 0$, and denote by $R_\varepsilon(t)$ the solution of the ODE $\frac{d}{dt} R_\varepsilon(t) = Q(R_\varepsilon(t)) + \varepsilon I$ with initial condition $R_\varepsilon(0) = R(0) + \varepsilon I$. The function $R_\varepsilon(t)$ is defined on some time interval $[0,T_\varepsilon)$. We claim that $R_\varepsilon(t) \in E$ for all $t \in [0,T_\varepsilon)$. To prove this, we argue by contradiction. Suppose that there exists a time $t \in [0,T_\varepsilon)$ such that $R_\varepsilon(t) \notin E$. Let 
\[\tau = \inf \{t \in [0,T_\varepsilon): R_\varepsilon(t) \notin E\}.\] 
Clearly, $\tau > 0$ and $R_\varepsilon(\tau) \in \partial E$. By Proposition \ref{algebraic.fact}, we can find an orthonormal four-frame $\{e_1,e_2,e_3,e_4\}$ and real numbers $\lambda,\mu \in [-1,1]$ such that 
\begin{align*} 
&R_\varepsilon(\tau)(e_1,e_3,e_1,e_3) + \lambda^2 \, R_\varepsilon(\tau)(e_1,e_4,e_1,e_4) \\ 
&+ \mu^2 \, R_\varepsilon(\tau)(e_2,e_3,e_2,e_3) + \lambda^2\mu^2 \, R_\varepsilon(\tau)(e_2,e_4,e_2,e_4) \\ 
&- 2\lambda\mu \, R_\varepsilon(\tau)(e_1,e_2,e_3,e_4) + (1 - \lambda^2) \, (1 - \mu^2) = 0. 
\end{align*} 
By definition of $\tau$, we have $R_\varepsilon(t) \in E$ for all $t \in [0,\tau)$. This implies 
\begin{align*} 
&R_\varepsilon(t)(e_1,e_3,e_1,e_3) + \lambda^2 \, R_\varepsilon(t)(e_1,e_4,e_1,e_4) \\ 
&+ \mu^2 \, R_\varepsilon(t)(e_2,e_3,e_2,e_3) + \lambda^2\mu^2 \, R_\varepsilon(t)(e_2,e_4,e_2,e_4) \\ 
&- 2\lambda\mu \, R_\varepsilon(t)(e_1,e_2,e_3,e_4) + (1 - \lambda^2) \, (1 - \mu^2) \geq 0 
\end{align*} 
for all $t \in [0,\tau)$. Hence, we obtain 
\begin{align*} 
&Q(R_\varepsilon(\tau))(e_1,e_3,e_1,e_3) + \lambda^2 \, Q(R_\varepsilon(\tau))(e_1,e_4,e_1,e_4) \\ 
&+ \mu^2 \, Q(R_\varepsilon(\tau))(e_2,e_3,e_2,e_3) + \lambda^2\mu^2 \, Q(R_\varepsilon(\tau))(e_2,e_4,e_2,e_4) \\ 
&- 2\lambda\mu \, Q(R_\varepsilon(\tau))(e_1,e_2,e_3,e_4) + \varepsilon \, (1 + \lambda^2) \, (1 + \mu^2) \leq 0. 
\end{align*} 
On the other hand, since $R_\varepsilon(\tau) \in E$, we have 
\begin{align*} 
&Q(R_\varepsilon(\tau))(e_1,e_3,e_1,e_3) + \lambda^2 \, Q(R_\varepsilon(\tau))(e_1,e_4,e_1,e_4) \\ 
&+ \mu^2 \, Q(R_\varepsilon(\tau))(e_2,e_3,e_2,e_3) + \lambda^2\mu^2 \, Q(R_\varepsilon(\tau))(e_2,e_4,e_2,e_4) \\ &- 2\lambda\mu \, Q(R_\varepsilon(\tau))(e_1,e_2,e_3,e_4) \geq 0 
\end{align*} 
by Proposition \ref{invariance.1}. This is a contradiction. 

Thus, we conclude that $R_\varepsilon(t) \in E$ for all $t \in [0,T_\varepsilon)$. It follows from standard ODE theory that $T \leq \liminf_{\varepsilon \to 0} T_\varepsilon$ and $R(t) = \lim_{\varepsilon \to 0} R_\varepsilon(t)$ for all $t \in [0,T)$. Therefore, we have $R(t) \in E$ for all $t \in [0,T)$. This completes the proof. \\

As in \cite{Bohm-Wilking}, we define a family of linear transformations $\ell_{a,b}$ on the space of algebraic curvature operators by 
\[\ell_{a,b}(R) = R + b \: \text{\rm Ric}_0 \owedge \text{\rm id} + \frac{a}{n} \, \text{\rm scal} \: \text{\rm id} \owedge \text{\rm id}.\] 
Here, $\text{\rm scal}$ and $\text{\rm Ric}_0$ denote the scalar curvature and trace-free Ricci tensor of $R$, respectively. Moreover, $\owedge$ denotes the Kulkarni-Nomizu product, i.e. 
\[(A \owedge B)_{ijkl} = A_{ik} \, B_{jl} - A_{il} \, B_{jk} - A_{jk} \, B_{il} + A_{jl} \, B_{ik}.\]
Using a result of C.~B\"ohm and B.~Wilking \cite{Bohm-Wilking}, we obtain: 

\begin{proposition} 
\label{invariance.3}
Assume that $b \in \big ( 0,\frac{\sqrt{2n(n-2)+4}-2}{n(n-2)} \big ]$ and $2a = 2b + (n-2) b^2$. 
Then the set $\ell_{a,b}(E)$ is invariant under the ODE $\frac{d}{dt} R = Q(R)$.
\end{proposition} 

\textbf{Proof.} 
By work of B\"ohm and Wilking (cf. \cite{Bohm-Wilking}, Theorem 2), it suffices to show that the set $E$ is invariant under the ODE $\frac{d}{dt} R = Q(R) + D_{a,b}(R)$, where $D_{a,b}(R)$ is defined by 
\begin{align*} 
D_{a,b}(R) &= ((n-2) \, b^2 - 2(a-b)) \, \text{\rm Ric}_0 \owedge \text{\rm Ric}_0 \\ 
&+ 2a \, \text{\rm Ric} \owedge \text{\rm Ric} + 2b^2 \, \text{\rm Ric}_0^2 \owedge \text{\rm id} \\ 
&+ \frac{nb^2 (1-2b) - 2(a-b)(1-2b+nb^2)}{n+2n(n-1)a} \, |\text{\rm Ric}_0|^2 \, \text{\rm id} \owedge \text{\rm id}.
\end{align*}
The first term on the right vanishes as $2a = 2b + (n-2) b^2$. By Corollary \ref{inclusions}, $E$ is a subset of $\tilde{C}$. Hence, every algebraic curvature operator $R \in E$ has nonnegative Ricci curvature. Consequently, we have $D_{a,b}(R) \geq 0$ for all $R \in E$. Since $E$ is invariant under the ODE $\frac{d}{dt} R = Q(R)$ by Proposition \ref{invariance.2}, we conclude that $E$ is also invariant under the ODE $\frac{d}{dt} R = Q(R) + D_{a,b}(R)$. \\

\section{Proof of the main theorem}

The proof of Theorem \ref{convergence.2} relies on the construction of a suitable pinching set. The concept of a pinching set was introduced in pioneering work of Hamilton (cf. \cite{Hamilton2}, Definition 5.1). B\"ohm and Wilking \cite{Bohm-Wilking} have a slightly more general notion of pinching set, which is more convenient for our purposes.

\begin{proposition}
\label{pinching.set}
Let $K$ be a compact set which is contained in the interior of $\tilde{C}$. 
Then there exists a closed, convex, $O(n)$-invariant set $F$ with the following properties: 
\begin{itemize} 
\item[(i)] $F$ is invariant under the ODE $\frac{d}{dt} R = Q(R)$. 
\item[(ii)] For each $\delta \in (0,1)$, the set $\{R \in F: \text{$R$ is not $\delta$-pinched}\}$ is bounded.
\item[(iii)] $K$ is a subset of $F$. 
\end{itemize}
\end{proposition}

\textbf{Proof.} 
By assumption, the set $K$ is contained in the interior of $\tilde{C}$. Using Proposition \ref{characterization.of.C.tilde}, we obtain 
\begin{align*} 
&R(e_1,e_3,e_1,e_3) + \lambda^2 \, R(e_1,e_4,e_1,e_4) \\ 
&+ \mu^2 \, R(e_2,e_3,e_2,e_3) + \lambda^2 \mu^2 \, R(e_2,e_4,e_2,e_4) - 2\lambda\mu \, R(e_1,e_2,e_3,e_4) > 0 
\end{align*} 
for all $R \in K$, all orthonormal four-frames $\{e_1,e_2,e_3,e_4\}$, and all pairs $(\lambda,\mu) \in 
\partial([-1,1] \times [-1,1])$. Hence, there exists a positive real number $N$ with the following properties: \\
1. We have 
\begin{align*} 
&R(e_1,e_3,e_1,e_3) + \lambda^2 \, R(e_1,e_4,e_1,e_4) \\ 
&+ \mu^2 \, R(e_2,e_3,e_2,e_3) + \lambda^2\mu^2 \, R(e_2,e_4,e_2,e_4) \\ 
&- 2\lambda\mu \, R(e_1,e_2,e_3,e_4) + N \, (1 - \lambda^2) \, (1 - \mu^2) > 0 
\end{align*} 
for all $R \in K$, all orthonormal four-frames $\{e_1,e_2,e_3,e_4\}$, and all pairs $(\lambda,\mu) \in [-1,1] \times [-1,1]$. \\
2. We have $\text{\rm tr}(R) \leq 2N$ for all $R \in K$. \\
Without loss of generality, we may assume that $N = 1$. Thus, $K$ is contained in the interior of the set $E$. Consequently, we can find real numbers $a \in (0,\frac{1}{2(n-1)}]$ and $b \in \big ( 0,\frac{\sqrt{2n(n-2) + 4} - 2}{n(n-2)} \big ]$ such that $2a = 2b + (n-2)b^2$ and $K \subset \ell_{a,b}(E)$. We now define $F_1 = \ell_{a,b}(E)$. Clearly, $F_1$ is closed, convex, and $O(n)$-invariant. Moreover, $F_1$ is invariant under the ODE $\frac{d}{dt} R = Q(R)$ by Proposition \ref{invariance.3}. 

We next consider the cones $\hat{C}(s)$ defined in \cite{Brendle-Schoen}. 
By continuity, we can find a real number $s_1 > 0$ such that $\ell_{a,b}(\hat{C}) \subset \hat{C}(s_1)$. Hence, it follows from Corollary \ref{inclusions} that 
\[\ell_{a,b}(R) + (1+2(n-1)a) \, I = \ell_{a,b}(R + I) \in \ell_{a,b}(\hat{C}) \subset \hat{C}(s_1)\] 
for all $R \in E$. Since $a \in (0,\frac{1}{2(n-1)}]$, we conclude that 
\[F_1 \subset \{R: R + 2 \, I \in \hat{C}(s_1)\}.\] 
Using Proposition 16 in \cite{Brendle-Schoen}, we can construct an 
increasing sequence of positive real numbers $s_j$, $j \in \mathbb{N}$, 
and a sequence of closed, convex, $O(n)$-invariant sets $F_j$, $j \in \mathbb{N}$, with the following properties: 
\begin{itemize}
\item[(a)] For each $j \in \mathbb{N}$, we have $F_{j+1} = F_j \cap \{R: R + 2^{j+1} I \in \hat{C}(s_{j+1})\}$. 
\item[(b)] For each $j \in \mathbb{N}$, we have $F_j \cap \{R: \text{\rm tr}(R) \leq 2^j\} \subset F_{j+1}$. 
\item[(c)] For each $j \in \mathbb{N}$, the set $F_j$ is invariant under the ODE $\frac{d}{dt} R = Q(R)$. 
\item[(d)] $s_j \to \infty$ as $j \to \infty$.
\end{itemize}
We now define $F = \bigcap_{j=1}^\infty F_j$. Clearly, $F$ is a closed, convex, $O(n)$-invariant set, which is invariant under the ODE $\frac{d}{dt} R = Q(R)$. Since $K \subset F_1 \cap \{R: \text{\rm tr}(R) \leq 2\}$, it follows from property (b) that $K \subset F_j$ for all $j \in \mathbb{N}$. Hence, $K$ is a subset of $F$. Finally, property (a) implies 
\[F \subset F_j \subset \{R: R + 2^j I \in \hat{C}(s_j)\}\] 
for all $j \in \mathbb{N}$. Since $s_j \to \infty$ as $j \to \infty$, the assertion follows from Proposition 15 in \cite{Brendle-Schoen}. \\

Having established the existence of a pinching set, the convergence of the normalized Ricci flow follows from work of Hamilton \cite{Hamilton2} (see also \cite{Bohm-Wilking}, Theorem 5.1):

\begin{theorem} 
\label{conv}
Let $(M,g_0)$ be a compact Riemannian manifold of dimension $n \geq 4$. Assume that the 
curvature tensor of $(M,g_0)$ lies in the interior of the cone $\tilde{C}$ for all points in $M$. 
Then the normalized Ricci flow with initial metric $g_0$ exists for all time and converges to a metric of constant sectional curvature as $t \to \infty$.
\end{theorem}

By Proposition \ref{characterization.of.C.tilde}, every curvature tensor satisfying (\ref{curvature.condition.2}) lies in the interior of the cone $\tilde{C}$. Thus, Theorem \ref{convergence.2} is an immediate consequence of Theorem \ref{conv}.

\end{document}